\documentclass{article}

\usepackage{amsfonts,amssymb,amsmath,amsthm}
\usepackage{fullpage}
\usepackage{color}
\usepackage{dsfont}
\usepackage{graphicx}
\usepackage{framed}
\usepackage[normalem]{ulem}
\usepackage{subfigure}
\usepackage{tikz}
\usetikzlibrary{decorations.pathreplacing}
\numberwithin{equation}{section}

\newtheorem{theorem}{Theorem}[section]

\newtheorem{remark}[theorem]{Remark}
\newtheorem{example}{Example}

\newcommand{\dx}{\, {\rm d} \mathbf x}

\newcommand{\ds}{\, {\rm d} s}
\newcommand{\bx}{\mathbf{x}}

\newcommand{\EE}{\mathbb{E}}
\newcommand{\VV}{\mathbb{V}}

\begin{document}

\title{Improved Multilevel Monte Carlo Methods for Finite Volume Discretisations of Darcy Flow in Randomly Layered Media}
\author{M. Park\footnote{School of Mathematical Sciences, University of Nottingham, University Park, Nottingham, NG7 2RD, UK. \tt{min.park@nottingham.ac.uk}}  \; and A.L. Teckentrup\footnote{Mathematics Institute, Zeeman Building, University of Warwick, Coventry, CV4 7AL, UK. \tt{a.teckentrup@warwick.ac.uk}}}
\date{}

\maketitle

\begin{abstract} We consider the application of multilevel Monte Carlo methods to steady state Darcy flow in a random porous medium, described mathematically by elliptic partial differential equations with random coefficients. The levels in the multilevel estimator are defined by finite volume discretisations of the governing equations with different mesh parameters. To simulate different layers in the subsurface, the permeability is modelled as a piecewise constant or piecewise spatially correlated random field, including the possibility of piecewise log-normal random fields. The location of the layers is assumed unknown, and modelled by a random process. We prove new convergence results of the spatial discretisation error required to quantify the mean square error of the multilevel estimator, and provide an optimal implementation of the method based on algebraic multigrid methods and a novel variance reduction technique termed Coarse Grid Variates.  
\end{abstract}

\section{Introduction}

Mathematical models of physical processes are frequently used for simulation. The parameters appearing in these models are often subject to uncertainty, due to for example lack of measurements or data, and the result of the simulation hence also becomes uncertain. In this work, we are interested in the simulation of a steady state groundwater flow, governed by Darcy's law, in a porous medium of which the permeability is not fully known. Assigning a suitable probability distribution to the permeability, the goal of the simulations will be to compute moments of quantities of interest related to the resulting pressure head and Darcy flux.

In recent years, multilevel Monte Carlo (MLMC) methods have frequently been applied and analysed in the context of partial differential equations with random coefficients \cite{cgst11,tsgu13,teckentrup12,bsz11,gkss11,anst14,canst14,ant14,ehm14}. Originally introduced by Heinrich \cite{heinrich01} in the context of parameter dependent integral equations and by Giles \cite{giles08} in the context of stochastic differential equations arising in mathematical finance, multilevel Monte Carlo methods present a significant computational saving compared to standard Monte Carlo methods. Exploiting the linearity of expectation, multilevel Monte Carlo methods combine a large number of samples of the quantity of interest at low spatial resolution with a small number of samples at finer spatial resolutions, to provide an accurate estimate of the quantity of interest at fine spatial resolution at low computational cost.

The purpose of this paper is to extend the theoretical convergence results available for MLMC methods, as well as improve on their practical performance through a novel variance reduction technique, with a particular focus on aspects relevant to the groundwater flow application. On the theoretical side, we firstly provide new convergence results on the spatial discretisation error, required to bound the mean square error of the MLMC estimator, in the case of finite volume discretisations. In applications such as groundwater flow modelling, finite volume methods are often preferred over methods such as standard finite elements due to local mass conservation \cite{dagan}. Secondly, we extend the range of layered permeability models considered in \cite{tsgu13} by allowing also the location of the layers to be uncertain. For the sake of generality, throughout the analysis we do not assume uniform coercivity or boundedness of the permeability in terms of the random parameter, but instead follow the more general framework suggested in \cite{cst13,tsgu13}.

On the implementation side, we propose a novel variance reduction technique, termed Coarse Grid Variates, to further lower the computational cost of MLMC estimators. The technique is designed to work for stationary models of the permeability. Assuming that for each sample of the permeability, we compute its values at the cell centres of the finite volume grid, we argue that samples of the permeability on a coarse grid can be extracted from the information contained in a sample of the permeability on a finer grid. The exact number of coarse grid samples that can be extracted from a single fine grid sample depends on the ratio of degrees of freedom between the two grids. Using this observation, together with an averaging procedure, we formulate a new MLMC estimator which is provably unbiased and in numerical simulations shows a variance reduction of up to 2 orders of magnitude at less than twice the computational cost of the standard MLMC estimator. 

The outline of the remainder of this paper is as follows. In section \ref{sec:probset}, we present the mathematical model of interest, together with any assumptions that we make on its components. Section \ref{sec:fv_mlmc} starts with a brief overview of finite volume methods and MLMC estimators, before we in section \ref{ssec:fv_anal} provide a convergence analysis of MLMC estimators based on finite volume discretisations. The theoretical results are illustrated in two and three spatial dimensions in section \ref{ssec:fv_num}. Section \ref{sec:cgv} is devoted to the derivation and numerical simulation of the Coarse Grid Variates technique. Finally, section \ref{sec:conc} provides some conclusions.

\section{Problem setting}\label{sec:probset}
The study of groundwater flow is well established,
and there is general scientific consensus that in many situations Darcy's law can be expected to lead to an
accurate description of the flow \cite{demarsily,delhomme79,cgss00}. The classical equations governing 
(steady state) single phase subsurface flow consist of Darcy's law coupled with an 
incompressibility condition:
\begin{equation}\label{eq:mod}
\mathbf{q} + \mathbf K \nabla p = \mathbf{g} \quad \text{and} \quad 
\text{div} \; \mathbf{q} = 0, \quad
\text{in } \ D  \subset \mathbb{R}^d, \ d=1,2,3, 
\end{equation}  
subject to suitable boundary conditions. In physical terms, $p$ denotes the 
pressure head of the fluid, $\mathbf K$ is the permeability tensor, 
$\mathbf{q}$ is the filtration 
velocity (or Darcy flux) and $\mathbf{g}$ are the source terms. 

A typical approach to incorporating the uncertainty in $p$ and $\mathbf q$ is to model the permeability tensor as a random field $\mathbf K = \mathbf K(\mathbf x,\omega)$ on $D \times \Omega$, with $(\Omega,\mathcal F, \mathbb P)$ a probability space \cite{ddgrtv05,delhomme79}. The model \eqref{eq:mod} then becomes a system of partial differential equations (PDEs) with random coefficients, which can be written in second order form as 
\begin{equation}\label{eq:mod2}
- \text{div} \, (\mathbf K(\mathbf x, \omega) \nabla p(\bx, \omega)) = f(\bx), \qquad \text{in} \quad D,
\end{equation}
with $f = -\text{div} \; \mathbf g$. The solution $p$ itself is also a random field on $D \times \Omega$. For simplicity, we assume that the boundary conditions and the sources $\mathbf g$ are deterministic, and restrict ourselves to convex polygonal/polyhedral domains $D$. For the more general case, we refer the reader to \cite{tsgu13,teckentrup12}.

In this general form, solving \eqref{eq:mod2} is extremely challenging 
computationally, and in practice it is therefore common to use relatively simple 
models for $\mathbf K$ that capture the most important features of the subsurface geometry and are as faithful as possible to the measured data.  For simplicity, we restrict our attention to scalar valued models of the permeability, i.e. we replace the tensor $\mathbf K$ by a scalar valued function $k$. Tensor valued coefficients are considered in \cite{teckentrup12}. In this paper, we are particularly interested in discontinuous models of the permeability, which are of particular practical interest due to their ability to model the different layers present in the subsurface. To this end, we assume that the computational domain $D$ is partitioned into $m$ disjoint convex polygonal subdomains $\{D_i\}_{i=1}^m$, with the permeability $k$ continuous on each subdomain $D_i$, $i=1,...,m$. 

In applications the exact location of the different layers of the subsurface is often not known exactly, and hence also uncertain. We incorporate this additional uncertainty into our model by allowing the partitioning $\{D_i\}_{i=1}^m$ to also be random, independent of $k$, leading to a random partitioning $\{D_i(\omega)\}_{i=1}^m$, for $\omega \in \Omega$. As before, we assume that for each realisation of the random partitioning, i.e. for a fixed $\omega$, the $m$ subdomains are disjoint and convex polygonal. 

It remains to specify the model of the permeability employed in the subdomains. For each subdomain, denote the permeability $k$ restricted to $D_i$ by $k_i$. A model that has been studied extensively in subsurface flow applications is the log-normal distribution, which allows the permeability to vary over many orders of magnitude and guarantees that the permeability takes positive values for almost all realisations. We will in particular work with the following two models:

\begin{example}\label{ex:logn} {\em (Piecewise constant model)} In each subdomain $D_i$, we model the permeability as a log-normal random variable: $k_i(\mathbf x,\omega) = \exp(Z_i(\omega))$, where $Z_i \sim N(\mu_i, \sigma_i^2)$, $i=1, \dots, m$. The mean and variance of the random variable are allowed to vary between the subdomains. As a function of the spatial variable $\mathbf x$, each $k_i$ is constant.
\end{example}

\begin{example}\label{ex:logn_spat}{\em (Piecewise correlated field model)} In each subdomain $D_i$, we model the permeability as a log-normal spatially correlated random field: $k_i(\mathbf x,\omega) = \exp(g_i(\mathbf x, \omega))$, where $g_i$ is a stationary Gaussian random field with constant mean $\mu_i(\mathbf x) = \mu_i$ and an exponential two point covariance function:
\begin{equation}\label{eq:exp_cov}
C(x,y) = \EE[ (g_i(\mathbf x) - \mu_i) (g_i(\mathbf y) - \mu_i)] = \sigma_i^2 \exp(- \|\mathbf x-\mathbf y\|_r / \lambda_i).
\end{equation}
Here, $\| \cdot \|_r$ denotes the usual $l_r$ norm on $\mathbb R^d$, and typically $r=1,2$. The parameter $\lambda_i$ is known as the correlation length and $\sigma_i^2$ as the variance of the Gaussian field $g_i$. It can be shown that in this case $k_i$, as a function of the spatial variable $\mathbf x$, is H\"older continuous with exponent less than $1/2$, $k_i(\cdot, \omega) \in C^{t}(\overline D_i)$, for any $t < 1/2$.
\end{example}

Other models of the permeability are of course possible, and the results of this paper are readily applicable in a variety of situations. In particular, the analysis in this paper holds also in the case where the covariance function $C(x,y)$ in Example \ref{ex:logn_spat} is replaced by another member of the Mat\`ern class of covariances with smoothness parameter $\nu > 1/2$ \cite{cst13}.

\section{Finite volume methods and multilevel Monte Carlo sampling}\label{sec:fv_mlmc}
In this section, we will describe and analyse the numerical methods employed in this paper. We will start with a description of the finite volume method used for the spatial discretisation in section \ref{ssec:fv}, before we briefly recall the multilevel Monte Carlo method in section \ref{ssec:mlmc} and finally in section \ref{ssec:fv_anal} prove the convergence of the multilevel Monte Carlo method for the finite volume method considered in section \ref{ssec:fv}.

\subsection{Finite volume discretisation}\label{ssec:fv}


The starting point of finite volume discretisations is the second order formulation \eqref{eq:mod}. One then chooses a non-overlapping partitioning of the domain $D$ into boxes (or {\em volumes}) $\mathcal B_h$, where $h$ denotes the mesh width of the partition. Integrating equation \eqref{eq:mod} over each box $B \in \mathcal B_h$ leads to a set of algebraic equations
\begin{equation}\label{eq:def_fv}
- \int_{B} \mathrm{div}(k(\mathbf x, \omega) \nabla p(\mathbf x, \omega)) \dx = \int_B f(\mathbf x) \dx, \qquad \forall B \in \mathcal B_h.
\end{equation}
The volume integral on the left hand side is transformed into a boundary integral using the Divergence Theorem:
\begin{equation}\label{def:fv}
- \int_{\partial B} k(\mathbf x, \omega) \nabla  p(\mathbf x, \omega) \cdot \mathbf n  \ds = \int_B f(\mathbf x) \dx,
\end{equation}
where $n$ denotes the unit outward normal and $\partial B$ denotes the boundary of the box $B$.
The specific finite volume scheme is now determined by the choice of volumes $\mathcal B_h$, as well as how the integrals in \eqref{def:fv} are computed (exactly or by quadrature). 

We will in this paper consider cell-centred finite volume methods on uniform rectangular meshes. For illustrative purposes, let us describe this discretisation in more detail in the particular case where the computational domain is the two-dimensional unit square, $D=(0,1)^2$, and the boxes $\mathcal B_h$ are squares. The cases of rectangular domains, and one or three spatial dimensions, are treated analogously.

We start by subdividing $[0,1]^2$ uniformly into a mesh of $m \times m$ square cells and denote by $B_{i,j}$ the cell $\left( \frac{i-1}{m}, \frac{i}{m}\right) \times \left( \frac{j-1}{m}, \frac{j}{m}\right)$, $i,j=1, \dots, m$, and by $\mathbf x_{i,j}$ its centre. We have $\overline{\mathcal B}_h = \cup_{1 \leq i,j \leq m} \overline B_{i,j}$, with $m = 1/h$. Let $k_{i,j}$ and $f_{i,j}$ be the values of $k$ and $f$ at $\mathbf x_{i,j}$, respectively, and denote by $p_{i,j}^\mathrm{FV}$ our approximation to $p$ at $\mathbf x_{i,j}$. We will approximate the right hand side of \eqref{eq:def_fv} by the midpoint rule,
\[
\int_{B_{i,j}} f(\mathbf x) \dx \approx f_{i,j} / m^2.
\]
To approximate the left hand side of \eqref{eq:def_fv} we consider separately each edge of $\partial B_{i,j}$. The contribution from the edge between $B_{i,j}$ and $B_{i+1,j}$ is again approximated by the midpoint rule. A simple approximation to $k$ on the edge is  its value at the midpoint, $k_{i+\frac{1}{2},j}$. As approximation to the gradient $\nabla p \cdot n$ on the edge we use the central finite difference $(p_{i+1,j}^\mathrm{FV} - p_{i,j}^\mathrm{FV}) / h$. The contributions from the other edges are approximated similarly, leading to the following form of the $(i,j)$-th equation:
\begin{equation}\label{eq:fv_lineq}
- k_{i, j-\frac{1}{2}} \, p_{i,j-1}^\mathrm{FV} - k_{i-\frac{1}{2},j} \, p_{i-1,j}^\mathrm{FV} + \Sigma_{i,j} \,  p_{i,j}^\mathrm{FV} -  k_{i+\frac{1}{2},j}\,  p_{i+1,j}^\mathrm{FV}  -  k_{i, j+\frac{1}{2}}\,  p_{i, j+1}^\mathrm{FV} =  f_{i,j} / m^2,
\end{equation}
where $ \Sigma_{i,j} =  k_{i, j-\frac{1}{2}} + k_{i-\frac{1}{2},j} +  k_{i+\frac{1}{2},j} +  k_{i, j+\frac{1}{2}}$. A Neumann boundary condition, i.e. a prescribed flux $-k \nabla p \cdot n = \psi_N$, on any part of the outer boundary of $(0,1)^2$ is straightforward to incorporate. The respective flux term on the left hand side of \eqref{eq:fv_lineq} is simple replaced by $\psi_N$  evaluated at the midpoint of the edge. To enforce a Dirichlet boundary condition, i.e. a prescribed pressure $p = \psi_D$, we simple replace the central difference by a one-sided difference. 

An alternative approximation to \eqref{eq:fv_lineq}, which is often used in subsurface flow applications, is derived using the harmonic average $\bar k_{i+\frac{1}{2},j}$ of $k_{i,j}$ and $k_{i+1,j}$ as the value of $k$ on the edge between $B_{i,j}$ and $B_{i+1,j}$. The analysis in section \ref{ssec:fv_anal} can be applied also in this case, leading to the same convergence rates as in Theorem \ref{thm:fv_lp}.

The linear system of equations arising from the approximation \eqref{eq:fv_lineq} takes the standard five-point stencil form. It is sparse and highly ill-conditioned, but it can be solved  efficiently and robustly with algebraic multigrid (AMG) methods. In fact, we will see in section \ref{ssec:scalable} that an iterative solver based on an AMG preconditioned Conjugate Gradient (CG) method scales optimally even in 3 spatial dimensions.

\subsection{Multilevel Monte Carlo methods}\label{ssec:mlmc}
We now briefly review the ideas of the multilevel Monte Carlo (MLMC) technique. For more details, we refer the reader to \cite{giles07,cgst11}.

Suppose we are interested in finding the expected value of a functional $Q = \mathcal{G}(p)$, where $p$ is the solution to the Darcy flow equation \eqref{eq:mod2}. Examples of functionals $Q$ of interest include the value of the pressure head $p$ or the Darcy flux $-k \nabla p$ at a particular point in the computational domain $D$, or the outflow over parts of the boundary. Since $p$ can not be computed exactly, we in practice use a finite volume approximation of $Q$, denoted by $Q_h := \mathcal{G}(p_h^{FV})$.

The main idea behind the MLMC technique is now as follows. Consider simulations with different mesh widths $h_\ell$, chosen such that 
\begin{equation}\label{eq:refine_factor}
h_\ell = s^{-1} h_{\ell-1}, \qquad \mbox{for } \ell = 0,1,...,L,
\end{equation}
where $s$ is a positive integer.
In contrast to the standard Monte Carlo (MC) approach, which only uses samples of $Q_h$ generated on the finest level $L$, samples on all grid levels $\ell=0,\ldots,L$ are taken into account in MLMC to estimate statistical moments of solution. Using the linearity of expectation, and denoting $Q_\ell := Q_{h_\ell}$, we have 
\begin{equation} \label{eq:telescopingsum}
\mathrm{E}[Q_{L}] = \mathrm{E}[Q_{0}] + \sum_{\ell=1}^L \mathrm{E}[Q_{\ell} - Q_{{\ell-1}}].
\end{equation}
Each of the expectations on the right hand side of \eqref{eq:telescopingsum} is now estimated independently using a Monte Carlo estimator, resulting in the MLMC estimator
\begin{equation}\label{eq:mlmc_est}
\widehat{Q}^{\mathrm{ML}}_L :=  \sum_{\ell = 0}^{L} \frac{1}{N_\ell} \sum_{i = 1}^{N_\ell} \left( Q_{\ell,i} - Q_{{\ell-1,i}}\right),
\end{equation}
where for simplicity we have set $Q_{-1} = 0$.
The number of Monte Carlo samples $N_\ell$ on each level is chosen such that the overall variance of the multilevel estimator is minimized for a fixed computational cost. It is important to note that the quantity $Q_{\ell,i} - Q_{{\ell-1,i}}$ in (\ref{eq:mlmc_est}) is computed from two discrete approximations with different mesh widths, but the same random sample $\omega^{(i)}$.

In order to quantify the accuracy of the multilevel estimator \eqref{eq:mlmc_est}, we consider the mean square error $\mathrm{MSE}(\widehat{Q},Q)$ of the estimator $\widehat{Q}$ as an estimator of $Q$:
\begin{equation*}
\mathrm{MSE}(\widehat{Q},Q) = \mathrm{V}[\widehat{Q}] + (\mathrm{E}[\widehat{Q}]-Q)^2.
\end{equation*}
Using the unbiasedness of Monte Carlo simulations, together with the fact that the $L+1$ individual Monte Carlo estimators in \eqref{eq:mlmc_est} are independent, the mean square error of the  MLMC estimator is
\begin{equation}\label{eq:mse_mlmc}
\mathrm{MSE}\left(\widehat{Q}^{\mathrm{ML}}_L, Q\right) = \sum_{\ell = 0}^{L}\frac{\mathrm{V}[Q_\ell - Q_{\ell-1}]}{N_\ell} + \left( \mathrm{E}[Q_L - Q]\right)^2.
\end{equation}


To achieve a mean square error $\mathrm{MSE}(\widehat{Q}^{\mathrm{ML}}_L,Q) $ at the tolerance level $\epsilon^2$, we evenly distribute $\epsilon^2$ between the two terms on the right hand side of (\ref{eq:mse_mlmc}). We furthermore denote the cost to compute one sample $Q_{\ell,i} - Q_{\ell-1,i}$ by $\mathcal C_\ell$, which includes both the cost of producing the sample of $k$ and solving the corresponding finite volume equations. We have the following results on the computational cost of the MC and MLMC estimators to achieve a mean square error of $\epsilon^2$.

\begin{theorem}\label{thm:mlmc_comp}
Suppose there exist positive constants $\alpha,\beta,\gamma, C_\alpha, C_\beta, C_\gamma >0$ such that $\alpha \geq \frac{1}{2}\min(\beta,\gamma)$ and 
\begin{equation*}
\begin{array}{ll}
(\mathit{M1})&|\mathrm{E}[Q_{\ell}-Q]| \leq C_\alpha h_\ell^{\alpha},\\
(\mathit{M2})&\mathrm{V}[Q_\ell - Q_{\ell-1}] \leq C_\beta h_\ell^{\beta},\\
(\mathit{M3})&\mathcal{C}_\ell \leq C_\gamma h_\ell^{-\gamma}.
\end{array}
\end{equation*}
Then for any $\epsilon < e^{-1}$, there exists a positive constant $C^{\mathrm{ML}}$, a value $L$ and a sequence $\{N_\ell\}_{\ell=0}^{L}$ such that $\mathrm{MSE}(\widehat{Q}_L^{\mathrm{\mathrm{ML}}},Q) < \epsilon^2$, and 
\begin{equation*} \label{eq:mlmc_cost}
\mathcal{C}_\epsilon(\widehat{Q}_L^{\mathrm{\mathrm{ML}}}) = \left\{ 
\begin{array}{ll}
         C^{\mathrm{ML}}\epsilon^{-2}, & \mbox{if $\beta < \gamma$},\\
		C^{\mathrm{ML}}\epsilon^{-2} (\log \epsilon)^2, & \mbox{if $\beta = \gamma$},\\
        C^{\mathrm{ML}}\epsilon^{-2-(\gamma-\beta)/\alpha}, & \mbox{if $\beta > \gamma$},\end{array}
\right.
\end{equation*}
whereas 
\begin{equation*}
\mathcal{C}_\epsilon(\widehat{Q}_M^{\mathrm{\mathrm{MC}}}) = C^{\mathrm{MC}} \epsilon^{-2-\gamma/\alpha}
\end{equation*}
for some positive constant $C^{\mathrm{MC}}$.
\end{theorem}

A proof of the above Theorem can be found in \cite{cgst11}. Further reductions in the computational cost of the MLMC estimator may be possible by using an optimal, uneven splitting of the total error between the two error contributions in \eqref{eq:mse_mlmc} \cite{canst14}, and by using optimised mesh hierarchies $\{h_\ell \}_{\ell=0}^L$ \cite{anst14}. The rates $\alpha$, $\beta$ and $\gamma$ are application dependent. The rates $\alpha$ and $\beta$ generally depend on the spatial regularity properties of $Q$ and the numerical method used for the approximation $Q_\ell$. The rate $\gamma$, on the other hand, depends on the method of sampling and the method used to solve the linear system of equations for each sample. In the best case, we have $\gamma$ approximately equal to $d$, the spatial dimension of the problem.

\subsection{Finite volume error analysis}\label{ssec:fv_anal}
This section is devoted to proving assumptions M1 and M2 of the complexity theorem for finite volume discretisations of the model problem \eqref{eq:mod} with the permeability as described in Examples \ref{ex:logn} and \ref{ex:logn_spat}. This will be done by showing that the cell-centred finite volume discretisation described in section \ref{ssec:fv} in fact is equivalent to a finite element discretisation with a particular quadrature scheme used to assemble the stiffness matrix. The convergence of the finite volume discretisation then follows from the convergence of finite element discretisations proven in \cite{cst13,tsgu13,teckentrup12}. We would like to point out here that a wider range of finite volume methods can in fact be analysed by comparison to a related finite element method, see e.g. \cite{h89}.

For ease of presentation, we again consider the case $D =(0,1)^2$ in detail, and restrict our attention to homogeneous Dirichlet conditions $\psi_D=0$ on the entire boundary $\partial D$. A standard finite element approximation of model problem \eqref{eq:mod} starts with the weak formulation of \eqref{eq:mod}, obtained by multiplying the equation by a test function $v \in H^1_0(D)$, integrating over the computational domain $D$ and applying Green's formula: find $p(\cdot, \omega) \in H^1_0(D)$ such that
\begin{equation}\label{eq:def_weak}
\int_D k(x, \omega) \nabla p(\mathbf x, \omega) \cdot \nabla v(\mathbf x) \dx = \int_D f(\mathbf x) v(\mathbf x) \dx, \qquad \forall v \in H^1_0(D).
\end{equation}
Here, $H^1_0(D)$ is the usual Sobolev space of functions with square integrable weak derivatives that vanish on the boundary:
\[
H^1_0(D) := \{ v : \int_D |v|^2 + | \nabla v|^2 \dx < \infty \quad \text{and} \quad v|_{\partial D} = 0 \}.
\]
The finite element approximation to \eqref{eq:def_weak}, denoted by $p^{FE}_h$, is then defined by $p^{FE}_h(\cdot, \omega) \in V_h$ and 
\begin{equation}\label{eq:def_fe}
\int_D k(\mathbf x, \omega) \nabla p^\mathrm{FE}_h(\mathbf x, \omega) \cdot \nabla v_h(\mathbf x) \dx = \int_D f(\mathbf x) v_h(\mathbf x) \dx, \qquad \forall v_h \in V_{h},
\end{equation}
where $V_h$ is a suitably chosen finite dimensional subspace of $H^1_0(D)$. For our purposes, we shall choose $V_h$ to be a space of continuous, piecewise bilinear functions on $D$ that vanish on the boundary $\partial D$. To facilitate the comparison with the finite volume discretisation described in section \ref{ssec:fv}, we choose the mesh for the finite element discretisation such that the degrees of freedom in the finite element method coincide with the degrees of freedom of the finite volume discretisation, which are the centres $\mathbf x_{i,j}$ of the boxes $B_{i,j} = \left( \frac{i-1}{m}, \frac{i}{m}\right) \times \left( \frac{j-1}{m}, \frac{j}{m}\right)$. The degrees of freedom of the standard, piecewise bilinear finite element method are at the vertices of the mesh. This means that for a fixed $m=1/h$, the finite element mesh is given by $\overline{\tilde{\mathcal B}}_h  = \cup_{0 \leq i,j \leq m} \overline{\tilde B}_{i,j}$, where with the nodes $y^0 = 0, y^m = 1$ and $y^i = (i - 1/2)h$, the square elements are given by $\tilde B_{i,j} = (y^i, y^{i+1}) \times (y^j, y^{j+1})$.

Since we impose Dirichlet boundary conditions on $\partial D$, only the interior nodes of the finite element mesh are considered as degrees of freedom. Note that these nodes are located exactly at the cell centres $\mathbf x_{i,j}$, for $1 \leq i,j \leq m$. To solve equation \eqref{eq:def_fe}, we now choose a basis for the piecewise bilinear finite element space $V_h$. We will use the well-known hat functions $\phi_{i,j}$, whose support is contained in the elements neighbouring the node $\mathbf x_{i,j}$. In particular, for any point $(x_1, x_2) \in [0,1]^2$, we have
\[
\phi_{i,j}(x_1, x_2) = 
\begin{cases}
\frac{(x_1 - y^{i-1})(x_2 - y^{j-1})}{(y^i - y^{i-1})(y^j - y^{j-1})} & \text{if} \quad (x_1, x_2) \in \tilde B_{i-1,j-1},\\
\frac{(x_1 - y^{i-1})(y^{j} - x_2)}{(y^i - y^{i-1})(y^j - y^{j-1})} & \text{if} \quad (x_1, x_2) \in \tilde B_{i-1,j}, \\
\frac{(y^{i} - x_1)(x_2 - y^{j-1})}{(y^i - y^{i-1})(y^j - y^{j-1})} & \text{if}\quad  (x_1, x_2) \in \tilde B_{i,j-1}, \\
\frac{(y^{i} - x_1)(y^{j} - x_2)}{(y^i - y^{i-1})(y^j - y^{j-1})} & \text{if} \quad (x_1, x_2) \in \tilde B_{i,j}, \\
0 & \text{elsewhere}.
\end{cases}
\] 
The finite element solution $p^\mathrm{FE}_h$ can be expanded in the basis $\phi_{i,j}$ as $p^\mathrm{FE}_h = \sum_{1 \leq i,j \leq m} p^\mathrm{FE}_{i,j} \,  \phi_{i,j}(x)$, where the coefficients $p^\mathrm{FE}_{i,j}$ will depend the particular realisation, i.e. on $\omega$. For ease of presentation, we will from now on drop the dependence of $k$, $p$ and $p_h^\mathrm{FE}$ on $\omega$. By choosing $v_h = \phi_{k,l}$ in \eqref{eq:def_fe}, for $1 \leq k,l \leq m$, we obtain a linear system of equations for the unknown coefficients $p^\mathrm{FE}_{i,j}$:
\begin{equation}\label{eq:fe_lineq}
p^\mathrm{FE}_{i,j} \, \int_{D} k(x)  \nabla \phi_{i,j}(x) \cdot \nabla \phi_{k,l}(x) \dx = \int_D f(x) \phi_{k,l}(x) \dx, \quad \text{for} \quad 1 \leq i,j,k,l \leq m.
\end{equation}
Due to the local support of the basis functions, the integrals appearing on the left hand side of \eqref{eq:fe_lineq} are zero if $|i - k| >1$ or $|j - l| > 1$. We will now devise a quadrature scheme to compute these integrals such that the linear equations \eqref{eq:fe_lineq} are the same as the linear equations \eqref{eq:fv_lineq} defining the finite volume approximation $p^\mathrm{FV}_{h}$. We will denote the corresponding solution of the finite element equations with quadrature by $p^\mathrm{qFE}_h$, with coefficients $p^\mathrm{qFE}_{i,j}$.

Let us start with the case $k=i+1, l=j$. Due to the local support of the basis functions, the integral over $D$ on the left hand side of \eqref{eq:fe_lineq} reduces to an integral over $\tilde B_{i,j-1} \cup \tilde B_{i,j}$. Expanding the dot product, this integral splits into the sum of two integrals. We approximate the first integral, which involves derivatives with respect to the first coordinate direction $x_1$, by the midpoint rule in $x_1$ and the trapezoidal rule in $x_2$:
\[
\int_{\tilde B_{i,j}} k \frac{\partial \phi_{i,j}}{\partial x_1} \, \frac{\partial \phi_{i+1,j}}{\partial x_1}\ \dx \approx \frac{h^2}{2} \left(   k \frac{\partial \phi_{i,j}}{\partial x_1} \, \frac{\partial \phi_{i+1,j}}{\partial x_1} ( x_1^{(i+1/2)}, x_2^{(j)} ) +  k \frac{\partial \phi_{i,j}}{\partial x_1} \, \frac{\partial \phi_{i+1,j}}{\partial x_1} ( x_1^{(i+1/2)}, x_2^{(j+1)} )    \right).
\]
The integral over $\tilde B_{i,j-1}$ is approximated in the same way, with the integrand evaluated at the two points $ ( x_1^{(i+1/2)}, x_2^{(j-1)} )$ and $ ( x_1^{(i+1/2)}, x_2^{(j)} )$. The second integral, involving derivatives with respect to $x_2$, is similarly approximated by the midpoint rule in $x_2$ and the trapezoidal rule in $x_1$. Explicit computations now show that this leads to the same set of equations as \eqref{eq:fv_lineq} (see, for example \cite[Exercise 4.1.8]{ciarlet} and \cite[\S 3.3]{buckeridge_thesis}). 

The computations above are easily extended to general rectangular domains $D$ and to three (or one) spatial dimensions. The quadrature scheme which makes the approximate finite element and the finite volume solution equivalent, is the one which uses the midpoint rule in the coordinate direction in which the derivatives are taken, and the trapezoidal rule in the remaining coordinate directions.

We now have the following convergence result, which follows immediately from the above analysis, together with the convergence results for the error in approximating $p$ by $p_h^{qFE}$ proven in \cite{cst13,tsgu13}.

\begin{theorem}\label{thm:fv_lp} Let the permeability $k$ be as in Example \ref{ex:logn} or \ref{ex:logn_spat}, and suppose $f \in L^2(D)$. Then for any linear functional $\mathcal G$ bounded on $H^1_0(D)$, we have 
\[
\EE[(\mathcal G(p) - \mathcal G(p_h^{FV}))^q]^{1/q} \; \leq  \; C \, h^{s},
\] 
for any $s < 1/2$ and $q < \infty$.
\end{theorem}

The convergence rate $s$ in Theorem \ref{thm:fv_lp} depends on the spatial (Sobolev) regularity of the solution $p$, and is, in the case of discontinuous coefficients $k$ limited above by $1/2$. The regularity results proved in \cite{tsgu13} only considered coefficients that are piecewise continuous with respect to a deterministic partitioning $\{D_i\}_{i=1}^m$ of $D$, but the same arguments can be applied to conclude on the regularity of $p$ in the case of a random partitioning. If we assume that for each realisation, the $m$ subdomains are disjoint and convex polygonal, and, in the case of zero Dirichlet conditions, no more than 2 subdomains meet at the boundary $\partial D$ and no more than three subdomains meet in the interior, we obtain the optimal convergence rate $s < 1/2$ in Theorem \ref{thm:fv_lp}.

From Theorem \ref{thm:fv_lp}, it now immediately follows from the triangle and reverse triangle inequalities, together with $\VV[X] \leq \EE[X^2]$ for any random variable $X$, that assumptions $\mathit{(M1)}$ and $\mathit{(M2)}$ are satisfied with $\alpha < 1/2$ and $\beta < 1$, respectively.

\begin{remark}\label{rem:quad} (Convergence rates for exact finite element solution) {\em In the presence of quadrature, the rate in Theorem \ref{thm:fv_lp} is optimal. However, the convergence rate of the exact finite element error $\EE[(\mathcal G(p) - \mathcal G(p_h^{FE}))^q]^{1/q}$ is twice that of Theorem \ref{thm:fv_lp}, with the error in functionals $\mathcal G$ converging with rate $2s$, for any $s < 1/2$. This faster convergence rate additionally requires $\mathcal G$ to be bounded on $L^2(D)$. We will see in section \ref{ssec:fv_num} that we in practice often observe this faster convergence rate even with quadrature.}
\end{remark}

\begin{remark}\label{rem:nonlin} (Nonlinear functionals) {\em The convergence result in Theorem \ref{thm:fv_lp} can be shown to hold also for non-linear functionals that are Fr\'echet differentiable, under suitable assumptions on the Fr\'echet derivative. For more details, see \cite{tsgu13}.}
\end{remark}

\begin{remark} (Boundary conditions) {\em The incorporation of non-zero Dirichlet and Neumann conditions is handled differently in finite element and finite volume methods. A result of the form in Theorem \ref{thm:fv_lp} for general boundary conditions can again be obtained, by applying the analysis above in interior boxes $\tilde B_{i,j}$, together with an argument bounding the error over the boundary boxes $\tilde B_{i,j}$ (see for example \cite[Section 3.2]{cst13}).}
\end{remark}


\subsection{Numerical examples}\label{ssec:fv_num}
In the following two sections, we present numerical results to verify the analysis provided in earlier sections, and to show that in fact all three assumptions of Theorem \ref{thm:mlmc_comp} are satisfied for the model problem considered in this paper. In the finite volume discretisation, we use the harmonic average of the permeability values in two neighbouring cells to approximate the permeability on the boundary.

\begin{figure}[ht!]
\centering
\includegraphics[width=0.45\textwidth]{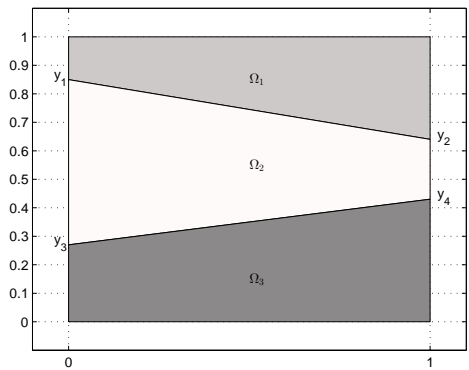}
\caption{A realisation of the random layers with four uniform random variables $y_1 \sim U(0.8,0.9), y_2 \sim U(0.6,0.7), y_3 \sim U(0.2,0.3)$ and $y_4 \sim U(0.4,0.5)$.}
\label{fig:random_layer}
\end{figure}

We consider a simple model of a randomly layered medium in two and three spatial dimensions, characterised by three different layers. First consider the two dimensional case, which is illustrated in Figure \ref{fig:random_layer}. To sample a realisation of the random medium for a given $\omega$, we draw four uniform random variables $y_1 \sim U(0.8,0.9), y_2 \sim U(0.6,0.7), y_3 \sim U(0.2,0.3)$ and $y_4 \sim U(0.4,0.5)$, and then draw straight lines between the points $(0, y_1)$ and $(1,y_2)$, and $(0,y_3)$ and $(1,y_4)$, respectively. Note that, as shown in Figure \ref{fig:random_layer}, for our particular choice of parametrisation, the subdomains created this way are always convex. For 3D simulations, we extrude 2D straight lines along the $x_3$-axis to create a 3D layered medium.

We consider two different model problems. \\

\textbf{Model Problem 1 (point evaluation of pressure):} Firstly, we consider a point value of the pressure. Since in higher dimensional space with $d = 2, 3$, point evaluation is not a bounded functional on $H^1(D)$, we regularise this type of functional by approximating the point value by a local average,
\begin{equation*}
M^{(1)}(v) := \frac{1}{|D^*|}\int_{D^*} v \dx
\end{equation*}
where $D^*$ is a small subdomain of $D$ that contains $x^*$ (\cite{gs02}). As governing equation, we consider the PDE
\begin{equation*}\label{eq:pde1}
-\nabla \cdot (k(\mathbf x,\omega)\nabla p(\mathbf x,\omega)) = 1, \quad \textnormal{ for } \mathbf x \in (0,1)^d 
\end{equation*}
with homogeneous Dirichlet boundary conditions. \\

\textbf{Model Problem 2 (outflow through boundary):} Secondly, we consider the PDE  
\begin{equation}\label{eq:pde2}
\nabla \cdot (k(\mathbf x,\omega)\nabla p(\mathbf x,\omega)) = 0, \quad \textnormal{ for } \mathbf x \in (0,1)^d
\end{equation}
with mixed boundary conditions $p\vert_{x_1 = 0} = 1, p\vert_{x_1 = 1} = 0$ and zero Neumann conditions on the remainder of the boundary.
We take as quantity of interest the outflow through the boundary $x_1 = 1$
\begin{equation*}\label{eq:qoi2}
M^{(2)}(v) =  -\int_{\partial D}\phi(\mathbf x) k(\mathbf x,\omega)\nabla v(\mathbf x,\omega) \cdot \mathbf n \ds
\end{equation*}
where $\phi$ is a weighting function such that $\phi\vert_{\{x_1=1\}} = 1$ and $\phi\vert_{\{x_1 = 0\} } = 0$. Notice that $M^{(2)}(p)$ is indeed equal to the outflow over the boundary $x_1 = 1$ because of the Neumann boundary conditions imposed on $p$.

To estimate the errors, we approximate the exact solution $u$ by a reference solution $u_{h^*}$ on a grid with mesh size $h^* = 1/512$ in 2D and $h^* = 1/128$ in 3D.

\begin{figure}
	\centering
	\subfigure{
	\includegraphics[width=.4\textwidth]{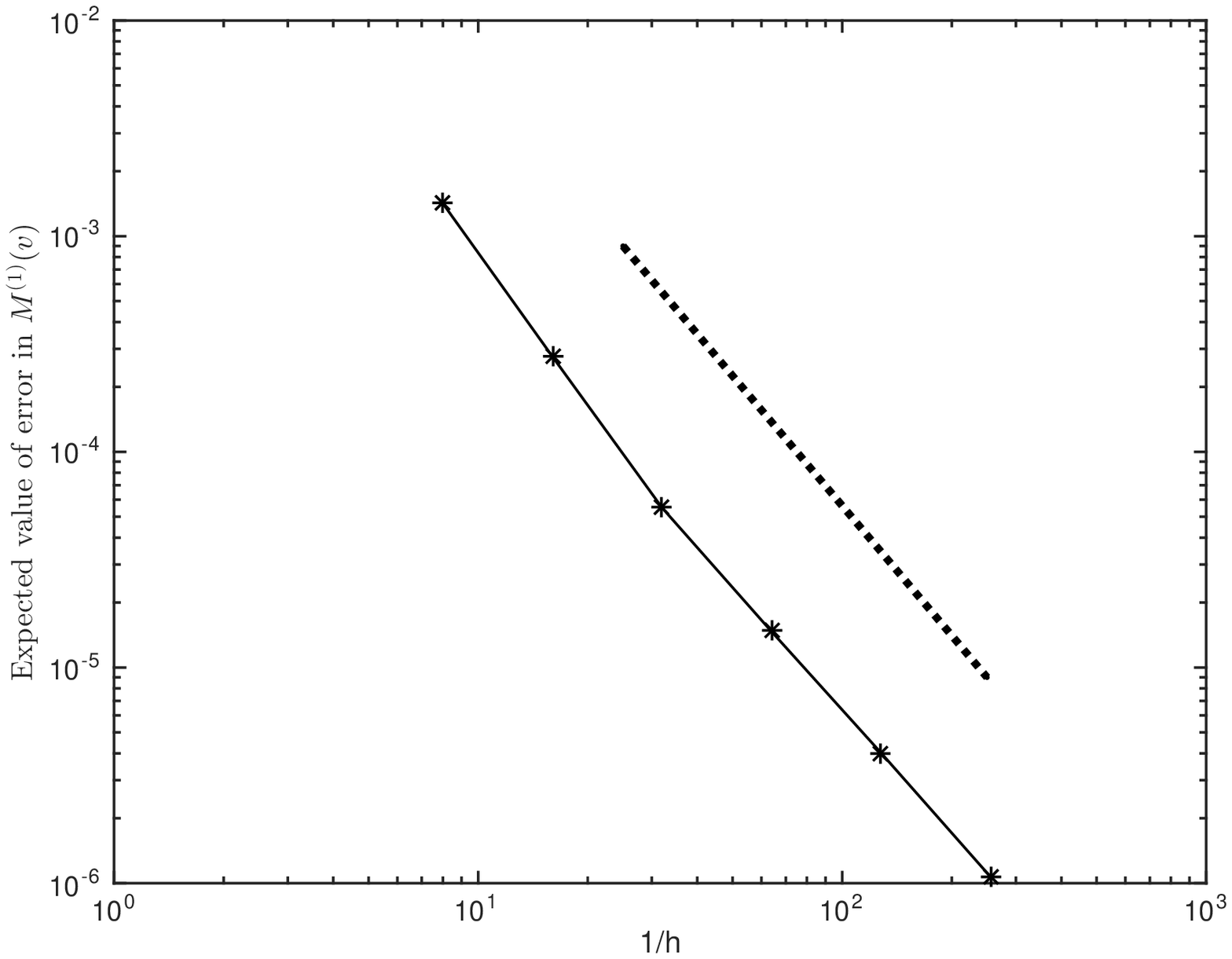}
	\includegraphics[width=.4\textwidth]{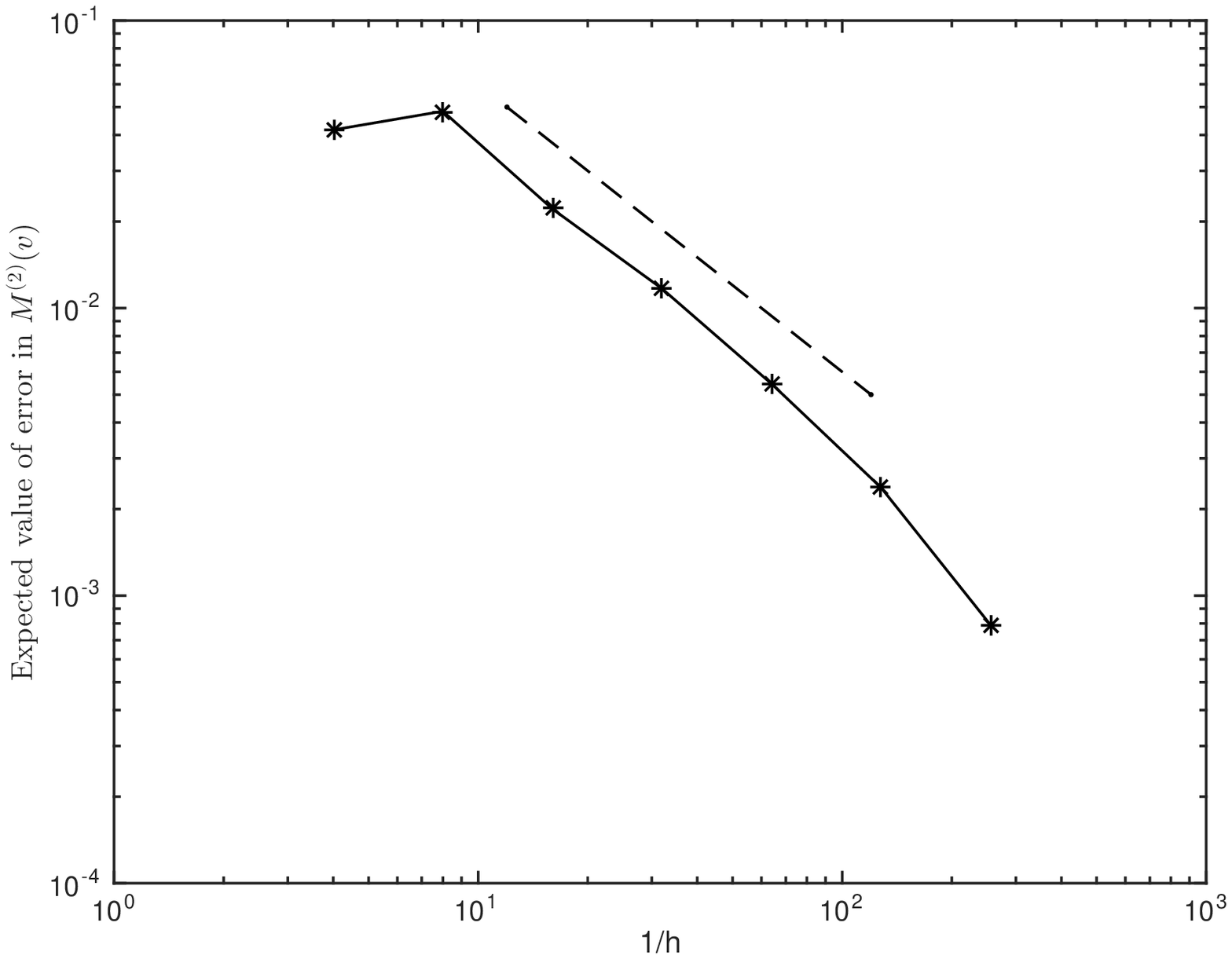}	
	}
	\subfigure{
	\includegraphics[width=.4\textwidth]{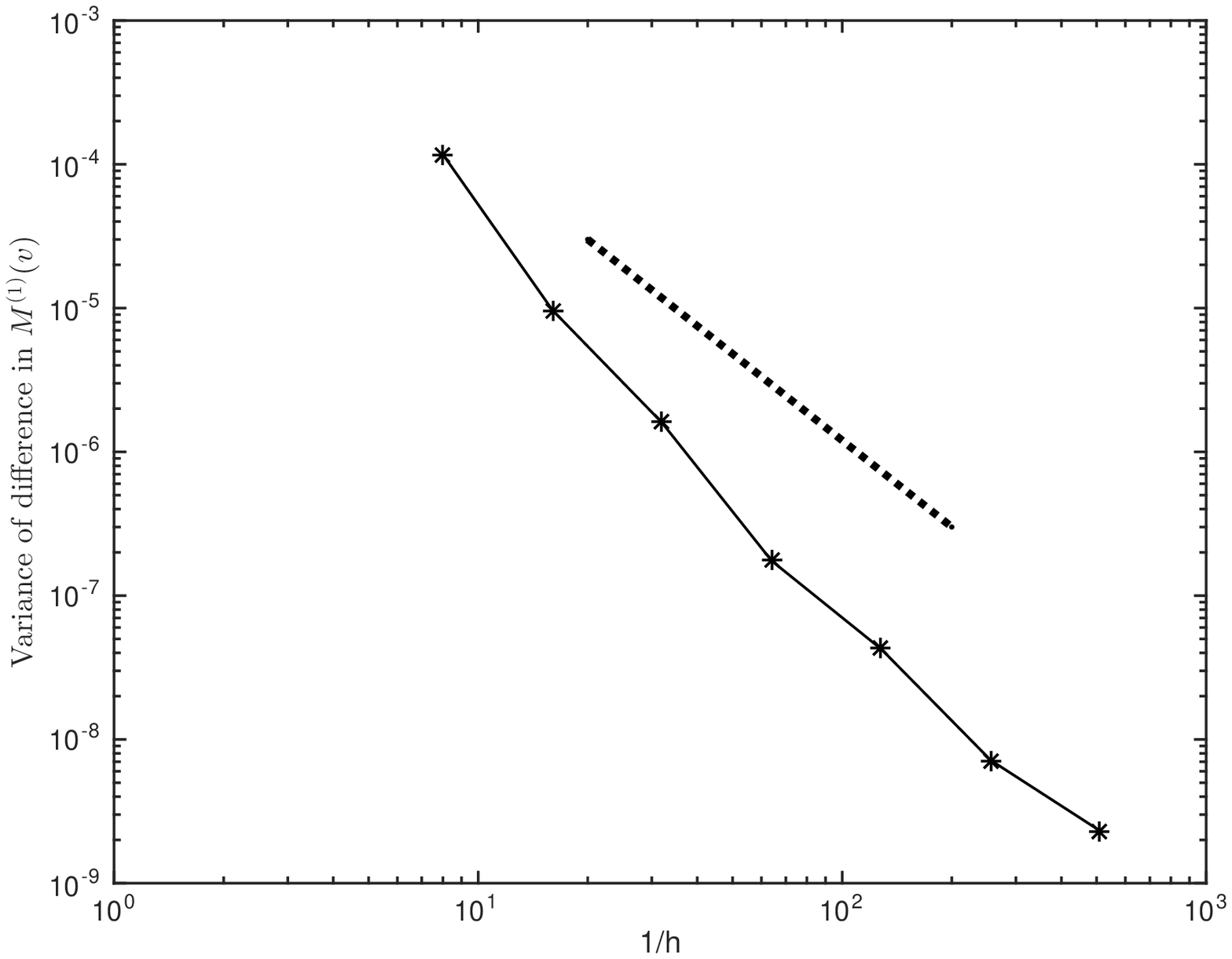}
	\includegraphics[width=.4\textwidth]{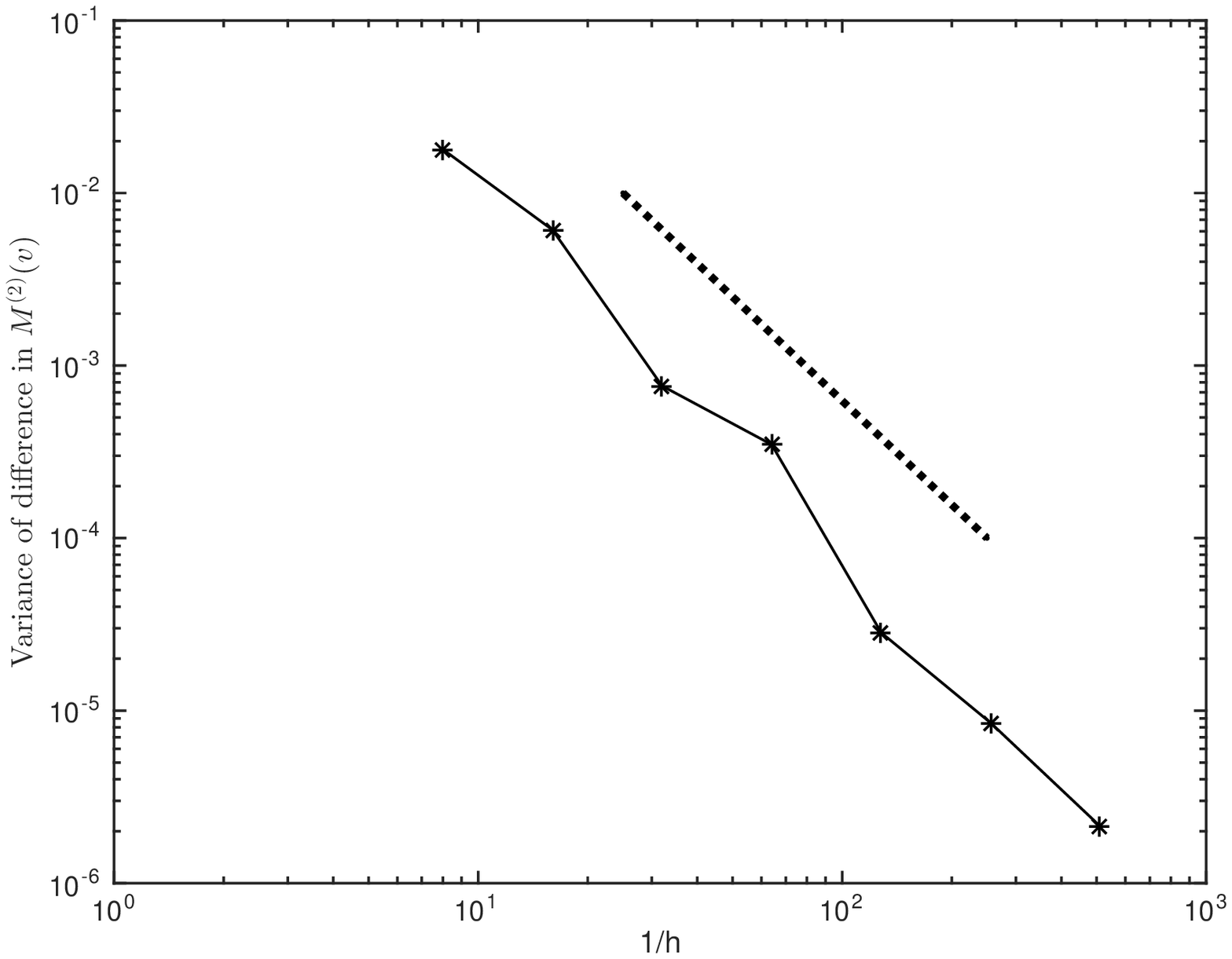}
	}	
	\caption{Results for piecewise constant permeability on random subdomains in 2D. Left-top plot: $|\EE[M^{(1)}(u_{h^*})-M^{(1)}(u_h)]|$ versus $1/h$ for Model Problem 1. Right-top plot: $|\EE[M^{(2)}(u_{h^*})-M^{(2)}(u_h)]|$ for Model Problem 2. Left-bottom plot: $\VV[M^{(1)}(u_h)-M^{(1)}(u_{2h})]$ versus $1/h$ for Model Problem 1. Right-bottom plot: $\VV[M^{(2)}(u_h)-M^{(2)}(u_{2h})]$ versus $1/h$ for Model Problem 2. The gradient of the dotted (resp. dashed) line is -2 (resp. -1).}
	\label{fig:const_random_layer}	
\end{figure}

\begin{figure}
	\centering
	\subfigure{
	\includegraphics[width=.4\textwidth]{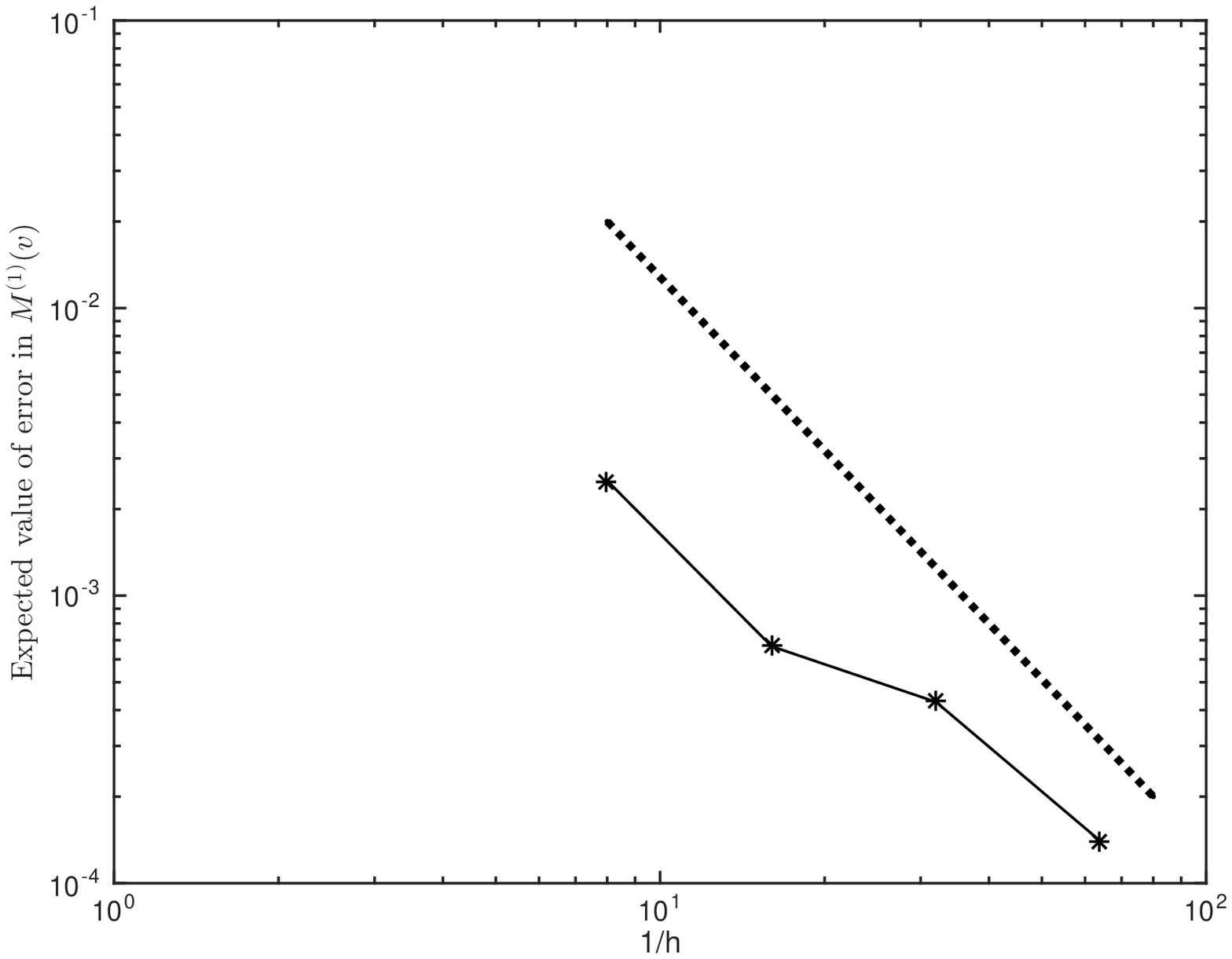}
	\includegraphics[width=.4\textwidth]{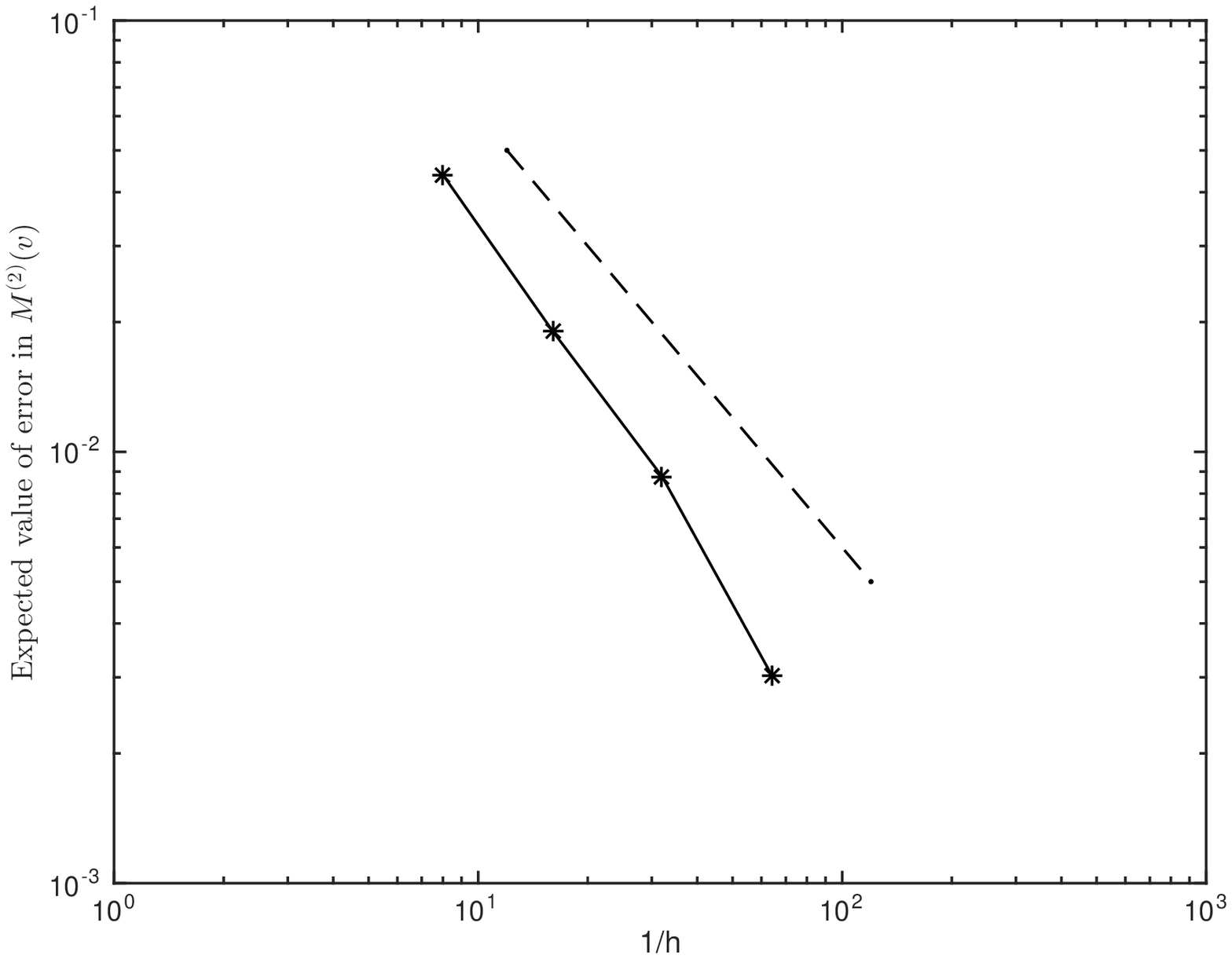}
	}
	\subfigure{
	\includegraphics[width=.4\textwidth]{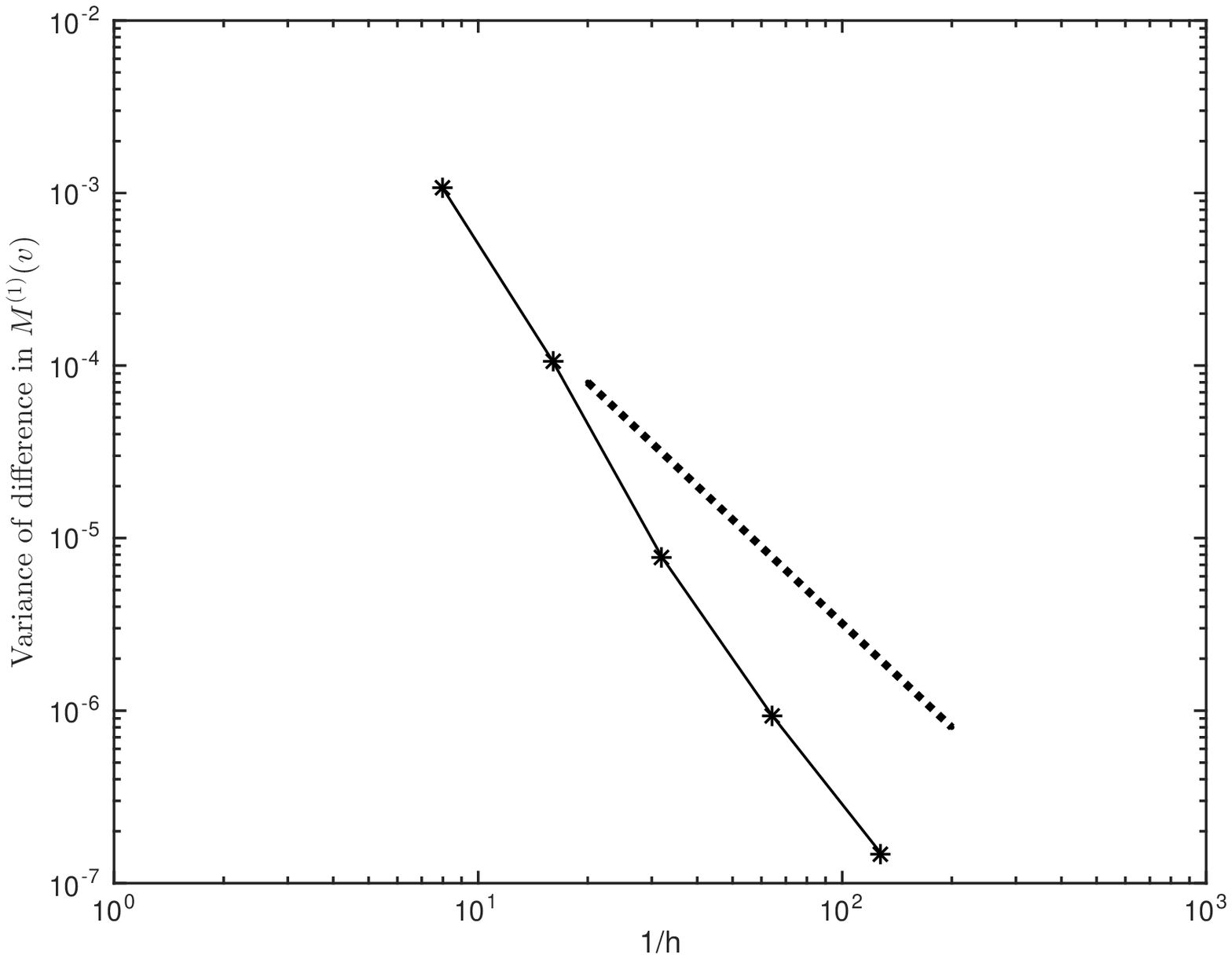}
	\includegraphics[width=.4\textwidth]{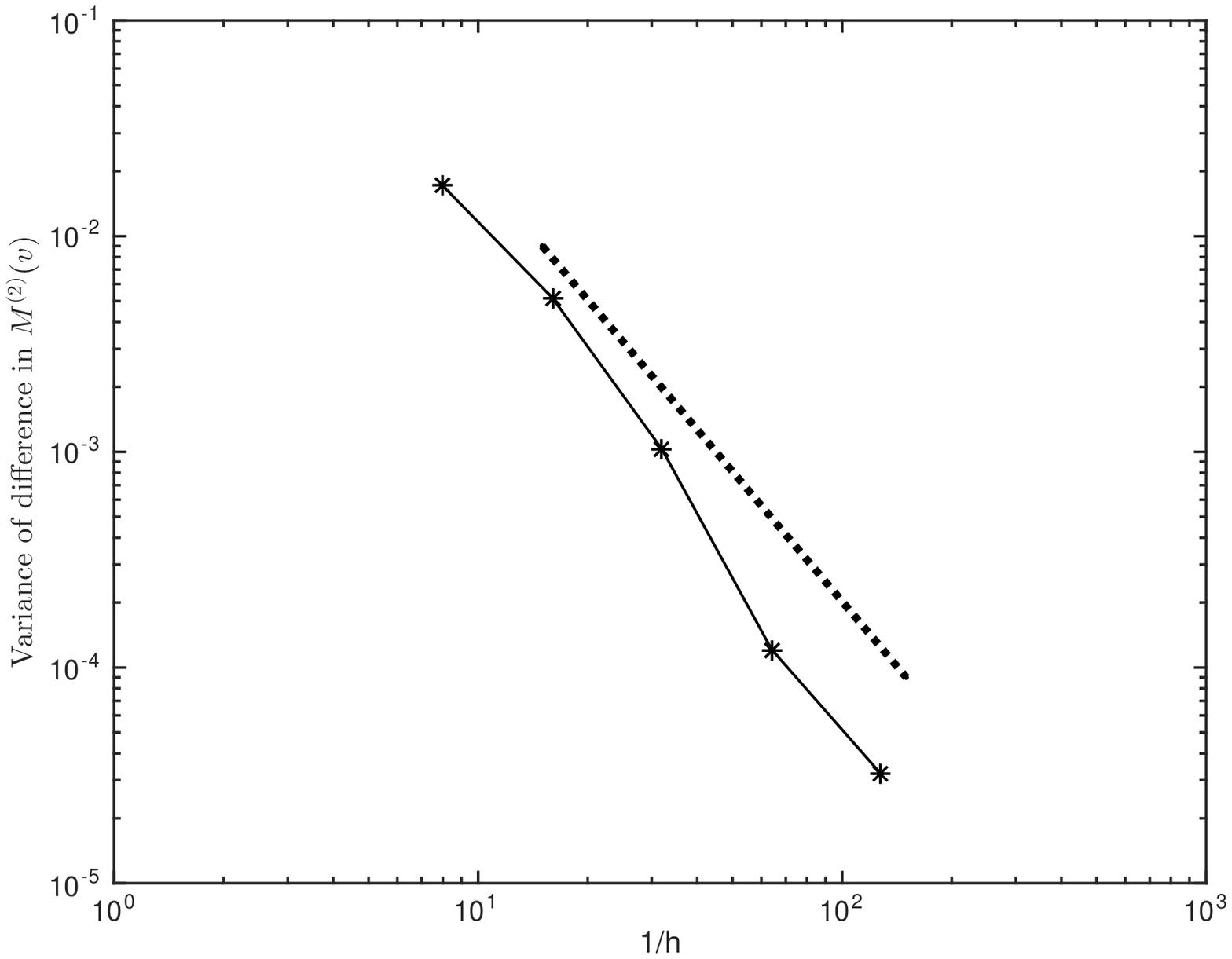}
	}	
	\caption{Results for piecewise constant permeability on random subdomains in 3D. Left-top plot: $|\EE[M^{(1)}(u_{h^*})-M^{(1)}(u_h)]|$ versus $1/h$ for Model Problem 1. Right-top plot: $|\EE[M^{(2)}(u_{h^*})-M^{(2)}(u_h)]|$ for Model Problem 2. Left-bottom plot: $\VV[M^{(1)}(u_h)-M^{(1)}(u_{2h})]$ versus $1/h$ for Model Problem 1. Right-bottom plot: $\VV[M^{(2)}(u_h)-M^{(2)}(u_{2h})]$ versus $1/h$ for Model Problem 2. The gradient of the dotted (resp. dashed) line is -2 (resp. -1).}
	\label{fig:const_random_layer3d}	
\end{figure}

Figures \ref{fig:const_random_layer} and  \ref{fig:const_random_layer3d} show results for the piecewise constant permeability model in Example 1. As described previously, we divide the medium into three (random) horizontal layers. For each sample of the layers, we then sample from three independent, standard normal random variables $z_1, z_2$ and $z_3$, and set the permeability values in the three regions to $\exp(z_1), \exp(z_2)$ and $\exp(z_3)$, respectively. 

In Figure \ref{fig:const_random_layer}, the two plots in the top row show quadratic convergence in $h$ for the error in expected value of $M^{(1)}$, and linear convergence for $M^{(2)}$. Theorem \ref{thm:fv_lp} and Remark \ref{rem:quad} suggest square-root convergence in the presence of quadrature, and linear convergence in the absence of quadrature, and so these convergence rates are faster than expected. In the two bottom plots, we observe faster than quadratic convergence for $\VV[M^{(1)}(u_h)-M^{(1)}(u_{2h})]$, and quadratic convergence for $\VV[M^{(2)}(u_h)-M^{(2)}(u_{2h})]$. In Figure \ref{fig:const_random_layer3d}, we observe the same convergence rates in three spatial dimensions. 

\begin{figure}
	\centering
	\subfigure{
	\includegraphics[width=.48\textwidth]{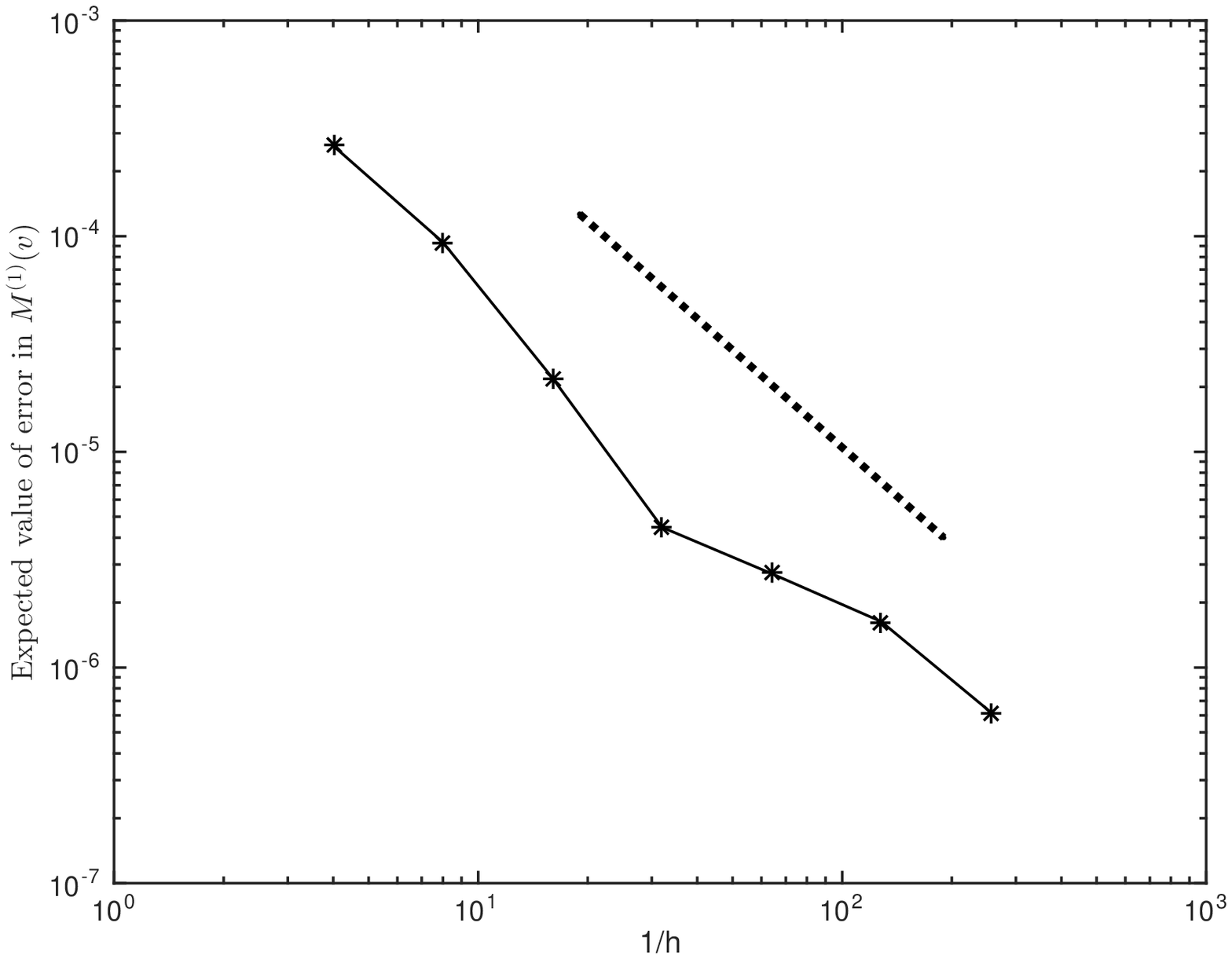}
	}
	\subfigure{
	\includegraphics[width=.48\textwidth]{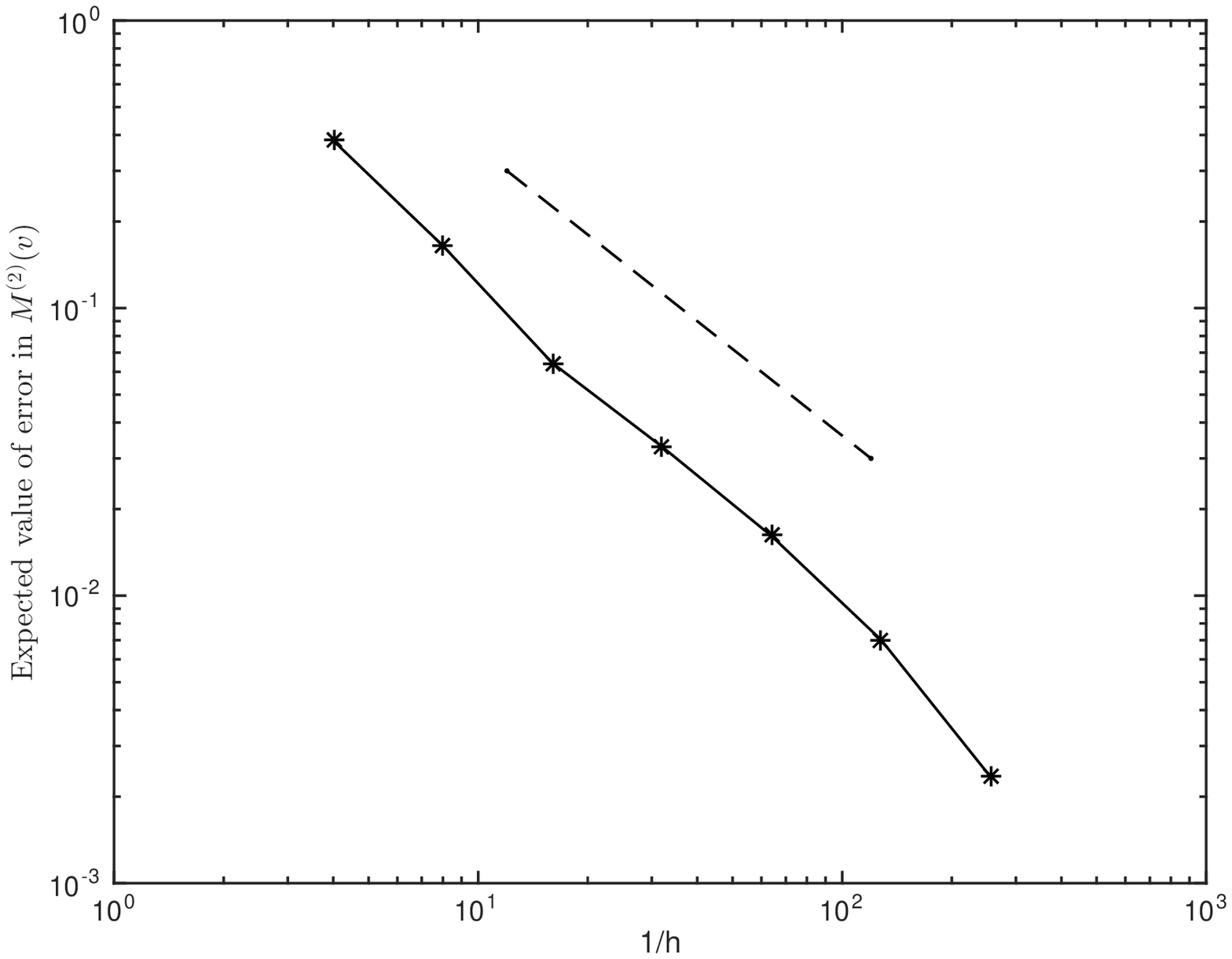}
	}	
	\caption{Results for piecewise correlated field permeability on random subdomains in 2D. Left plot: $|\EE[M^{(1)}(u_{h^*})-M^{(1)}(u_h)]|$ versus $1/h$ for Model Problem 1. Right plot: $|\EE[M^{(2)}(u_{h^*})-M^{(2)}(u_h)]|$ for Model Problem 2.  The gradient of the dotted (resp. dashed) line is -3/2 (resp. -1).}
	\label{fig:cov_random_layer}	
\end{figure}

Figure \ref{fig:cov_random_layer} shows results for the piecewise correlated field model from Example 2 in two spatial dimensions. As before, we divide $D = (0,1)^2$ into three (random) horizontal layers, and model the permeability in the 3 layers by 2 different log-normal distributions. The parameters in the top and bottom layer are taken to be $\mu_1 = 0, \lambda_1 = 0.3$ and $\sigma_1^2 = 1$, and for the middle layer we take $\mu_2 = 4, \lambda_2 = 0.1$ and $\sigma_2^2 = 1$, assuming no correlation across layers. We use the 2-norm exponential covariance function (\ref{eq:exp_cov}). To generate realisations of a log-normal random field, we use the circulant embedding technique, which is an exact and fast method of generating samples from stationary Gaussian random fields on a regular grid  \cite{wood:1994, dn97}. The complexity of this sampling method is $\mathcal O( h_\ell^{-d} \log h_\ell^{-d})$.
 
In Figure \ref{fig:cov_random_layer}, we observe $\mathcal{O}(h^{\frac{3}{2}})$ convergence for the error in expected value for $M^{(1)}$ and linear convergence for $M^{(2)}$. 

%

\subsection{Scalable linear solver ($\gamma = d$) in 3D}\label{ssec:scalable}
In this section, we investigate the scalability of the iterative linear solver, which determines the value of $\gamma$ in Theorem \ref{thm:mlmc_comp}. Since it is the most challenging, we only consider the case of three spatial dimensions. Results in one and two spatial dimensions are similar. As already mentioned in section \ref{ssec:mlmc}, the value of $\gamma$ also depends on the sampling method used. The Circulant Embedding method used in the previous section has log-linear complexity in the number of degrees of freedom, $\mathcal O( h_\ell^{-d} \log h_\ell^{-d})$. Out of the sampling and the linear solver, the linear solver usually is the more costly in terms of absolute CPU time.

The finite volume discretisation of a realisation of model problem \eqref{eq:mod2} results in a symmetric positive definite (SPD) stiffness matrix. As a linear solver for SPD systems, we use a preconditioned conjugate gradient (PCG) method, with the preconditioner given by one cycle of an algebraic multigrid (AMG) solver.

In 3D simulations, memory requirements for storing the stiffness matrices increases rapidly with the mesh width $h_\ell$, since the number of degrees of freedom in the finite volume method grows like $h_\ell^{-d}$. Thus, the efficiency of the linear solver strongly depends on the number of non-zeros in the matrices. In order to achieve optimal computational complexity, we use the aggressive coarsening technique described in \cite{trottenberg:2001}. The outer PCG iteration is run until the stopping criterion
\[
\|r_m\|/\|r_0\| < 10^{-10}
\]
is satisfied, where $r_m$ and $r_0$ are the residuals obtained at iterations $m$ and 0, respectively. 

\begin{figure}[ht!]
\centering
\includegraphics[width=0.6\textwidth]{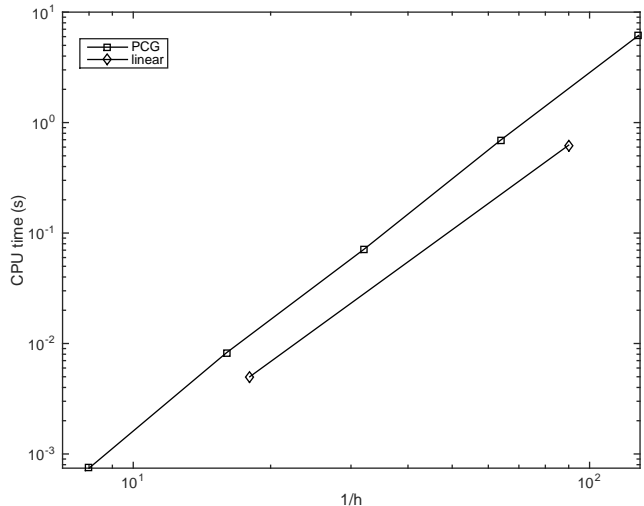}
\caption{Average CPU time (in seconds) for the solution of 1000 linear systems with random coefficients, solved by using AMG-PCG. }
\label{fig:scalability}
\end{figure}

Figure \ref{fig:scalability} shows the average CPU time (in seconds) for the solution of 1000 linear systems using AMG-PCG. The systems were generated by discretising the partial differential equation (\ref{eq:pde2}) as in Model Problem 2, with random coefficients given by a log-normal random field with 2-norm exponential covariance function \eqref{eq:exp_cov} with $\sigma^2 = 1$ and $\lambda = 0.3$. The CPU time increases linearly with the degrees of freedom. Hence, AMG-PCG is a scalable linear solver for 3D problems with random coefficients.

\section{Variance reduction through Coarse Grid Variates}\label{sec:cgv}
\begin{figure}
\centering
\begin{tikzpicture}
\draw[-][very thick](0,0)--(8,0);
\draw[very thick, dashed](8,0)--(12,0);
\draw[-][very thick](12,0)--(14,0);
\draw[thick] (0,-.1) node[below]{$x_0 = 0$} -- (0,0.1);
\draw[thick] (2,-.1) node[below]{$x_1$} -- (2,0.1);
\draw[thick] (4,-.1) node[below]{$x_2$} -- (4,0.1);
\draw[thick] (6,-.1) node[below]{$x_3$} -- (6,0.1);
\draw[thick] (8,-.1) node[below,xshift=.2cm]{$x_{4}$} -- (8,0.1);
\draw[thick] (12,-.1) node[below,xshift=.1cm]{$x_{m_\ell}$} -- (12,0.1);
\draw[thick] (14,-.1) node[below,xshift=.1cm]{$x_{m_\ell+1} = 1$} -- (14,0.1);
\draw [thick, black,decorate,decoration={brace,amplitude=10pt},xshift=0.4pt,yshift=-0.4pt](1.0,0.5) -- (5.0,0.5) node[black,midway,yshift=0.6cm] {\footnotesize $2h_\ell$};
\draw [thick, black,decorate,decoration={brace,amplitude=10pt,mirror},xshift=0.4pt,yshift=-0.4pt](3.0,-0.5) -- (7.0,-0.5) node[black,midway,yshift=-0.6cm] {\footnotesize $2h_\ell$};
\draw [thick, black,decorate,decoration={brace,amplitude=10pt,mirror},xshift=0.4pt,yshift=-0.4pt](0.0,-0.5) -- (2.0,-0.5) node[black,midway,yshift=-0.6cm] {\footnotesize $h_\ell$};
\node at (1.0,0)[above]{$k_\ell(1)$};
\node at (3.0,0)[below]{$k_\ell(2)$};
\node at (5.0,0)[above]{$k_\ell(3)$};
\node at (7.0,0)[below]{$k_\ell(4)$};
\node at (13.0,0)[above, xshift=.3cm]{$k_\ell(m_\ell)$};

\node at (1,0)[circle,fill=gray,inner sep=0pt,minimum size=4pt]{};
\node at (3,0)[circle,fill=gray,inner sep=0pt,minimum size=4pt]{};
\node at (5,0)[circle,fill=gray,inner sep=0pt,minimum size=4pt]{};
\node at (7,0)[circle,fill=gray,inner sep=0pt,minimum size=4pt]{};
\node at (13,0)[circle,fill=gray,inner sep=0pt,minimum size=4pt]{};
\end{tikzpicture}
\caption{1D uniform grid $\{x_0,\ldots,x_{m_\ell+1} \} \in \overline D = [0, 1]$, with $h_\ell = x_i - x_{i-1}$. Grey circles represent the locations of sampling points $k_\ell(i)$. } \label{fig:1D_grid}
\end{figure}
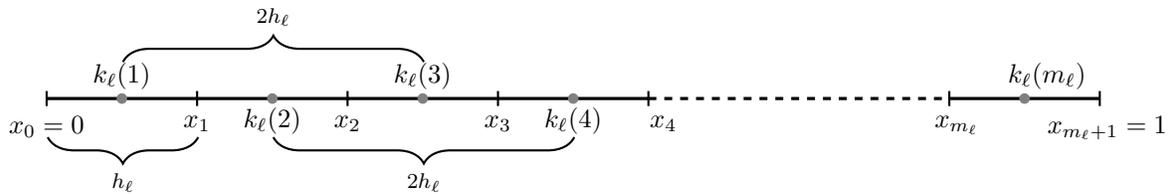

Different variance reduction techniques, such as antithetic variates, control variates and importance sampling, have been developed to increase the accuracy of and speed up the simulation process \cite{boyle:1997}. In this section, we introduce a new variance reduction technique specifically designed for multilevel Monte Carlo estimators, which we term Coarse Grid Variates (CGV). In particular, the variance reduction technique is designed to reduce the variance of the Monte Carlo estimators of the differences $Q_\ell - Q_{\ell-1}$, for $\ell=1, \dots, L$. We will describe the method in the context of finite volume discretisations of the Darcy flow equation as discussed in earlier sections, but it is easily adaptable to other applications.

\subsection{Description of method}
For ease of presentation, let us for the moment assume that the finite volume meshes used in the simulation are regular in the sense that nodes are uniformly distributed in each coordinate direction. We also assume that the number of cells in each coordinate direction is a power of 2, with $s$ = 2 in the growth condition (\ref{eq:refine_factor}). Under these assumptions, the hierarchy of finite volume meshes is nested. We will crucially also assume that the distribution of $k$ is {\em stationary} in the sense that it is invariant under shifts in $d$-dimensional space.

To compute a sample $Q_{\ell,i} - Q_{\ell-1,i}$, we need to produce a sample of the coefficient $k$ at the cell centres of the finite volume meshes on level $\ell$ and $\ell-1$. Let us denote these coarse and fine samples by $k_{\ell,i}$ and $k_{\ell-1,i}$, respectively. The aim is now to extract a sample $k_{\ell-1,i}$ from a given sample  $k_{\ell,i}$, using the fact that the distribution $k$ is stationary. 

We first illustrate our general idea in the one-dimensional setting. Figure \ref{fig:1D_grid} shows the computational domain $D=(0,1)$, divided into $m_\ell$ cells. The fine grid sample $k_{\ell,i}$ consists of $m_\ell$ values. The coarse grid sample $k_{\ell-1,i}$ should consist of $m_{\ell-1} = m_\ell/2$ values, and since the distribution of $k$ is stationary, it follows that the two subvectors of $k_{\ell,i}$ given by its even and odd entries, respectively, are both acceptable as samples $k_{\ell-1,i}$. With $k_{\ell}^{(o)} := [k_\ell(1)\;, k_\ell(3)\,; \cdots \;, k_\ell(m_\ell-1)]$, and $k_{\ell}^{(e)} = [k_\ell(2)\;, k_\ell(4)\;, \cdots \;, k_\ell(m_\ell)]$,  both $Q_{\ell,i} - Q_{\ell-1,i}^{(o)}$ and $Q_{\ell,i} - Q_{\ell-1,i}^{(e)}$ are valid samples of $Q_{\ell} - Q_{\ell-1}$. Furthermore, we notice that  $Q_{\ell,i} - (Q_{\ell-1,i}^{(o)} +  Q_{\ell-1,i}^{(e)})/2$ is also a valid sample, and this is what we will use as a basis for our variance reduction technique.

In general $d$ dimensional space, the random field $k_\ell$ on level $\ell$ has $2^d$ sub-vectors $k_\ell^{(j)}$, for $j = 1, \ldots ,2^d$, which have the distribution required for the random field $k_{\ell-1}$ on level $\ell-1$. Given a sample $k_{\ell,i}$ of $k_\ell$, we then define the averaged coarse grid quantity
\begin{equation}\label{eq:cgv_est_ddim}
Q_{\ell-1}^{\mathrm{cgv}} = \frac{1}{2^d}\sum_{j=1}^{2^d} {Q}(k_{\ell}^{(j)}),
\end{equation}
and finally the estimator
\begin{equation}\label{eq:mc_cgv_est}
\widehat{Y}^{\mathrm{cgv}}_{\ell,N_\ell} := \frac{1}{N_\ell} \sum_{i = 1}^{N_\ell} \left( Q_{\ell,i} - Q_{{\ell-1,i}}^{\mathrm{cgv}}\right)
\end{equation}
for $Y_\ell = Q_\ell - Q_{\ell-1}$. Note that this is still an unbiased estimator of $\EE[Y_\ell]$. As we will see in section \ref{ssec:cgv_num}, the estimator \eqref{eq:mc_cgv_est} has a considerably smaller variance than the standard MC estimators for $Q_\ell - Q_{\ell-1}$ based on only one sample $Q_{\ell-1,i}$. This is due to the fact that the CGV estimator \eqref{eq:mc_cgv_est} incorporates much more information due to the averaging of samples on the coarse level.

The computational cost of the multilevel estimator based on a telescoping sum of estimators \eqref{eq:mc_cgv_est} will be less than twice that of of the standard MLMC estimator presented in section \ref{ssec:mlmc}. For each difference $Q_\ell - Q_{\ell-1}$, the total cost of the additional computations required is no more than the cost of a single solve on the fine level $\ell$, since we have $\gamma \geq d$. In practice, the sampling cost of the Coarse Grid Variates is less than that of 2 fine-grid samples because no random field needs to be generated on the coarse grid. Instead, the random vector already generated on the fine grid is used on the coarse grid.

%

\subsection{Numerical examples}\label{ssec:cgv_num}
In this section we provide numerical results demonstrating the effectiveness of the CGV method. We consider the mixed boundary value problem Model Problem 2 already considered in section \ref{ssec:fv_num}, with quantity of interest the outflow functional $M^{(2)}$. We model $k$ as a continuous log-normal random field, such that $\log k$ has 1-norm exponential covariance function (\ref{eq:exp_cov}) with $\lambda = 0.3$ and $\sigma^2 = 1$. 

\begin{figure}
	\centering
	\includegraphics[width=.6\textwidth]{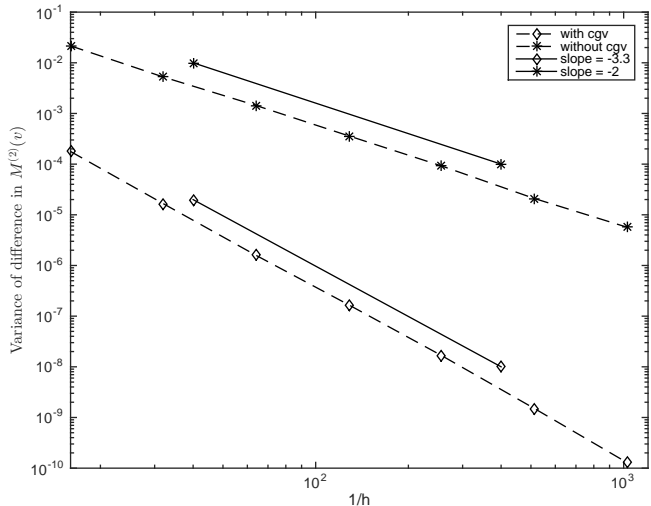}
	\caption{Variance plot of $Y_\ell$ for Model Problem 2 in 2D,with continuous log-normal permeability.}
	\label{fig:2d_cgv}	
\end{figure}
Figure \ref{fig:2d_cgv} shows results for the model problem in two spatial dimensions. In particular, it shows the variance of $Y_\ell$ with and without the use of Coarse Grid Variates. The dashed line with diamonds indicates the results using CGV. We see that the use of CGV leads to a significantly lower variance, both in terms of the constant and even a faster decay rate with mesh width $h_\ell$. Already on the coarsest grid of width $h_\ell = 1/8$, the variance of $Y_\ell$ is reduced by a factor of almost 200. These results suggest that the CGV can further lower the rate of growth of the computational cost of the MLMC algorithm.
\begin{figure}
	\centering
	\includegraphics[width=.6\textwidth]{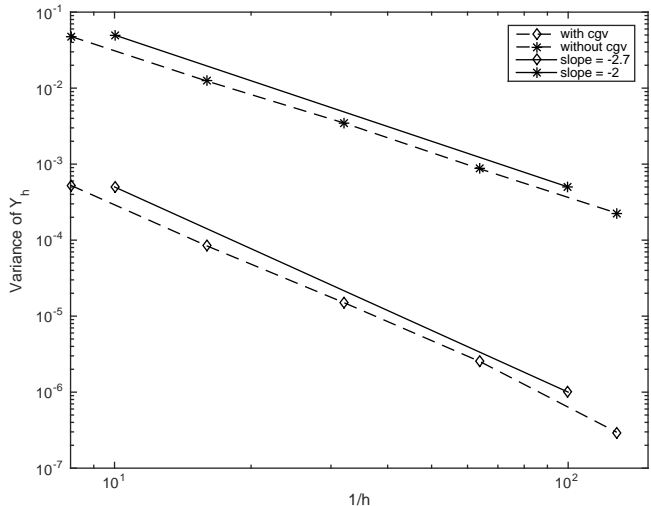}	
	\caption{Variance plot of $Y_\ell$ for Model Problem 2 in 3D,with continuous log-normal permeability.}
	\label{fig:3d_cgv}	
\end{figure}

In Figure \ref{fig:3d_cgv}, we observe a very similar behaviour of $Y_\ell$ for the three dimensional model problem. On the coarsest grid $h_\ell = 1/8$, the variance of $Y_\ell$ is reduced by a factor of 100 when using CGV. We again also observe a faster decay rate of $Y_\ell$ with $h_\ell$ when using CGV.

In numerical simulations with other quantities of interest, such as the point evaluation considered in Model Problem 1 from section \ref{ssec:fv_num}, the benefits of using the CGV approach were again clear, although not as pronounced as in Figures \ref{fig:2d_cgv} and \ref{fig:3d_cgv}. The variance reduction attainable from the CGV approach seems to depend on the specific model problem and quantity of interest. However, in our tests, the CGV approach always resulted in gains over standard MLMC.

\section{Conclusions}\label{sec:conc}
In this work, we considered the application of multilevel Monte Carlo methods to elliptic PDEs with log-normal random coefficients. We extended the existing theory to cover simple finite volume discretisations of the governing equations, and model problems where the location of different layers in the subsurface is also subject to uncertainty. Theoretical results were confirmed by numerical simulations. Finally, we proposed a new variance reduction technique, termed Coarse Grid Variates, designed for the multilevel Monte Carlo method. Numerical simulations show a variance reduction around 2 orders of magnitude compared to standard multilevel Monte Carlo for model problems in two and three spatial dimensions.

\section*{Acknowledgements}
The first author was partially supported by the EPSRC grant EP/K031430/1. The second author was partially supported by EPSRC grant EP/K034154/1.

\bibliographystyle{siam}
\bibliography{bibMLMC}
\end{document}